\documentstyle[doublespace,11pt]{article}

\font\sm=cmr10

   \csname @addtoreset\endcsname{equation}{section}
    \csname @addtoreset\endcsname{figure}{section}
    \csname @addtoreset\endcsname{table}{section}

\newcounter{thanksnum}
\def\thanksnumber#1
{\setcounter{thanksnum}{\value{footnote}}\setcounter{footnote}{#1}%
                     \addtocounter{footnote}{-1}\footnotemark
                     \setcounter{footnote}{\value{thanksnum}}}
\def\newtheoremz#1{\@ifnextchar[{\@othmz{#1}}{\@nthmz{#1}}}

\def\@nthmz#1#2{%
\@ifnextchar[{\@xnthmz{#1}{#2}}{\@ynthmz{#1}{#2}}}

\def\@xnthmz#1#2[#3]{\expandafter\@ifdefinable\csname #1\endcsname
{\@definecounter{#1}\@addtoreset{#1}{#3}%
\expandafter\xdef\csname the#1\endcsname{\expandafter\noexpand
  \csname the#3\endcsname \@thmcountersepz \@thmcounterz{#1}}%
\global\@namedef{#1}{\@thmz{#1}{#2}}\global\@namedef{end#1}{\@endtheoremz}}}

\def\@ynthmz#1#2{\expandafter\@ifdefinable\csname #1\endcsname
{\@definecounter{#1}%
\expandafter\xdef\csname the#1\endcsname{\@thmcounterz{#1}}%
\global\@namedef{#1}{\@thm{#1}{#2}}\global\@namedef{end#1}{\@endtheoremz}}}

\def\@othmz#1[#2]#3{\expandafter\@ifdefinable\csname #1\endcsname
  {\global\@namedef{the#1}{\@nameuse{the#2}}%
\global\@namedef{#1}{\@thmz{#2}{#3}}%
\global\@namedef{end#1}{\@endtheoremz}}}

\def\@thmz#1#2{\refstepcounter
    {#1}\@ifnextchar[{\@ythmz{#1}{#2}}{\@xthmz{#1}{#2}}}

\def\@xthmz#1#2{\@begintheoremz{#2}{\csname the#1\endcsname}\ignorespaces}
\def\@ythmz#1#2[#3]{\@opargbegintheoremz{#2}{\csname
       the#1\endcsname}{#3}\ignorespaces}

\def\@thmcounterz#1{\noexpand\arabic{#1}}
\def\@thmcountersepz{.}
\def\@begintheoremz#1#2{ \trivlist \item[\hskip \labelsep{\bf #1\ #2}]}
\def\@opargbegintheoremz#1#2#3{ \trivlist
      \item[\hskip \labelsep{\bf #1\ #2\ (#3)}]}
\def\@endtheoremz{\endtrivlist}

\newtheorem{theorem}{Theorem}[section]

\newtheorem{corollary}{Corollary}[section]
\newtheorem{definition}{Definition}[section]

\newtheorem{proposition}{Proposition}[section]

\def\abstract{\if@twocolumn
\section*{\abstractname}%
\else \small
\quotation
\fi}

\def\endabstract{\if@twocolumn\else\endquotation\fi}

\def\defi{\stackrel{{\scriptscriptstyle \Delta}}{=}}

\def\w{\widehat}

\def\R{{\bf R}}
\def\E{{\bf E}}
\def\P{{\bf P}}

\def\s{\delta}

\def\s{\sigma}

\def\p{\partial}

\newcommand{\be}{\begin{equation}}
\newcommand{\ee}{\end{equation}}
\newcommand{\bd}{\begin{displaymath}}
\newcommand{\ed}{\end{displaymath}}
\newcommand{\ba}{\bd\begin{array}{rl}}
\newcommand{\ea}{\end{array}\ed}
\font\sm=cmr10
\date{ March, 1999  }
\def\BS{{\scriptscriptstyle BS}}
\title {A paradox of diffusion market model related with existence of
winning combinations of options\thanks{ Supported by Russian
Foundation for Basic Research grant 99-01-00886}}
\author{Nikolai Dokuchaev\thanks{Correspondence address:
The Institute of Mathematics and Mechanics, St.Petersburg State
University, Bibliotechnaya pl.2, Petrodvoretz, St.Petersburg,
198904,  Russia. E-mail: dokuchaev@pobox.spbu.ru; fax
7(812)4286998}\\ {\sm The Institute of Mathematics and Mechanics},
{\sm St.Petersburg State University, Russia}}
\begin{document} \maketitle
\begin{abstract}
We consider strategies of investments into options and diffusion
market model. It is shown that there exists a correct proportion
between "put" and "call" in the portfolio such that the average
gain is almost always positive for a generic Black and Scholes
model. This gain is zero if and only if the market price of risk
is zero. It is discussed a paradox related to the corresponding
loss of option's seller.
\end{abstract}
\section{Introduction}
There are many works devoted investment strategies dealing with
risky assets (stocks). We consider strategies of investments into
options for a generic stochastic diffusion model of a financial
market. It is assumed that there is a risky stock and a risk-free
asset (bond), and that it is available European options "put" and
"call" on that stock at the initial time We consider only
strategies of selecting options portfolio at the initial time. The
selection of this portfolio id the only action of the investor;
after that, he/she waits until the expiration  time to accept gain
or loss.
\par
We show that there exist a correct proportion  between "put" and
"call" options in the portfolio such that the average gain is
almost always positive. This gain is zero if and only if the
market price of risk is zero, i.e. when the appreciation rate of
the stock is equal to the interest rate of the risk free asset
(i.e. $a\neq r$ in (\ref{(1)})--(\ref{(2)}) below). It is
discussed a paradox related to the corresponding loss of option's
seller.
\section{Definitions}
Consider the generic model of financial market   consisting of a
risk-free asset (bond, or bank account) with price $B_t$, and a
risky asset (stock) with price    $S_t$. We are given a standard
probability space with a probability measure $\P$ and a standard
Brownian motion $w_t$. The bond and stock prices evolve as
\be
\label{(1)}
B_t=e^{rt}B_0,
\ee
\be
dS_t=aS_t dt+\s S_tdw_t. \label{(2)} \ee Here $r\ge 0$ is the
risk-free interest rate, , $\s>0$ is the volatility,  and $a\in\R$
is the appreciation rate.  We assume that
$t\in [0,T]$,  where
$T>0$ is a given terminal time. The equation (\ref{(2)}) is
It\^o's equation, and can be rewritten as
$$ S_t=S_0\exp\left(at-\frac{\s^2t}{2}+\s w_t\right). $$
\par
We assume that it is available the options  "put"  and "call" on
that stock for that price defined by the Black-Scholes formula
(see e.g. Strong (1994)).
\par
Further, we assume that  $\s$, $r$, $B_0>0$ and $ S_0>0$ are
given, but the constant  $a$ is unknown.
\par Let $
p_{\BS}(S_0,K,r,T,\s)$ denote Black-Scholes price for "put"
option, and $ c_{\BS}(S_0,K,r,T,\s)$ denote Black-Scholes price
for "call" option. Here $S_0$ is the initial stock price,
 $K$ is the strike price,
$r$ is the risk-free interest rate,
 $\s$ is the volatility,
 $T$ is the expiration time.
\par
We remind the Black-Scholes formula.
Let
$$
\Phi(x)\defi\frac{1}{\sqrt{2\pi}}\int_{-\infty}^x e^{-\frac{y^2}{2}}dy,
$$
\be
d\defi
\frac{\ln\left(S_0/K\right)+T\left(r+\s^2/2\right)}{\s\sqrt{T}},\quad
d^-=d-\s\sqrt{T}.
\label{(3)}
\ee
Then
\be
c_{\BS}(S_0,K,r,T,\s)=S_0\Phi(d)-Ke^{-rT}\Phi(d^-),
\label{(4)}
\ee
\be
p_{\BS}(S_0,K,r,T,\s)=c_{\BS}(S_0,K,r,T,\s)-S_0+Ke^{-rT}
\label{(5)} \ee (see e.g. Strong (1994)-Duffie  (1988)).
\par
We denote by  $X_ 0$ the initial wealth of an investor
(i.e. at the initial $t=0$), and we denote by $X_T$
the wealth of the investor at the terminal time
 $t=T$.
\par
Consider a vector $\left(K_p,\mu_p,K_c,\mu_c\right)$, such that $K_p>0,\mu_p\ge
0,K_c>0,\mu_c\ge 0$. We shall use  this vector to describe
 the following strategy: buy a portfolio of options which
consists of  $\mu_p$ "put" options with the strike price $K_p$ and
$\mu_c$ options with the strike price $K_c$, with the same
expiration time  $T$; thus,
\be
X_0=\mu_pp_{\BS}(S_0,K_p,r,T,\s)+\mu_cc_{\BS}(S_0,K_p,r,T,\s),
\label{(6)}
\ee
(We assume that the options are available for the Black-Scholes
price). We have assumed that the investor does not take any
other actions until the expiration time. In that case,
the terminal wealth at time $t=T$ will be
\be
X_T=\mu_p(K_p-S_T)^+ +\mu_c(S_T-K_c)^+.
\label{(7)}
\ee
\begin{definition}
\label{def1}
The vector
$\left(K_p,\mu_p,K_c,\mu_c\right)$
is said to be a strategy.
\end{definition}
For a case of risk-free "hold-only-bond" strategy,
 $X_T=e^{rT}X_0$.
It is natural to compare the results of any investment with the
risk-free investment.
\begin{definition}
\label{def2}
The difference $\E X_T-e^{rT}X_0$ is said to be the average
gain.
\end{definition}
Note that the appreciation rate
$a$ in this definition is fixed but unknown.
The average gain for a strategy
depends on $a-r$. For example, for a "call" option holder, when
 ($\mu_p=0$), the average gain is positive
if  $a>r$.
\section{The result}
Let  $d_p$ and  $d_c$ be defined by
 (\ref{(3)}), where
$d$ is defined after substituting $K=K_p$ or $K=K_c$
correspondingly.
\begin{theorem}
\label{Th1}
Let $\mu_p>0$, $\mu_c>0$ and
\be
\frac{\mu_c}{\mu_p}=\frac{1-\Phi(d_p)}{\Phi(d_c)}.
\label{(8)}
\ee
Then the average gain for the strategy
 $\left(K_p,\mu_p,K_c,\mu_c\right)$
is positive for any $a\neq r$, i.e.
\be
\E X_T>e^{rT}X_0 \quad \forall a\neq r. \label{(9)} \ee Moreover,
\be
\E X_T=e^{rT}X_0\quad \hbox{ if}\quad a=r.
\label{(10)}
\ee
For any  $\left(K_p,K_c\right)$, the proportion  (\ref{(10)})
is the only one which ensures  (\ref{(8)}): for any other
proportion  $\mu_c/\mu_p$ there exists  $ \in \R$ such that the
average gain
is negative.
\end{theorem}
\begin{corollary}
\label{cor1} Let the variable $a$ be random, independent of
$w(\cdot)$ and such that $\P(a\neq r)>0$.  Then  $\E
X_T>e^{rT}X_0$ for the strategy  from Theorem \ref{Th1}.
\end{corollary}
\par
Set
$$
Q(x,t)\defi \mu_pp_{\BS}(x,K_p,r,T-t,\s)+\mu_cc_{\BS}(x,K_c,r,T-t,\s),
$$
\be\Delta (x,t)=\frac{dQ(x,t)}{dx}.
\label{(11)}
\ee
\begin{corollary}
\label{cor2} Let $\left(K_p,\mu_p,K_c,\mu_c\right)$ be such as in
Theorem \ref{Th1}, then
\be
\Delta(S_0,0)=0.
\label{(12)}
\ee
\end{corollary}
In terms of Strong (1994), the equality   (\ref{(12)}) means that
the position is risk-neutral  at time $t=0$.
\par
{\bf Example.} Consider example with same parameters as in Strong
(1994)[p. 109]. Let $K_p=K_c=\$25$, $S_0=\$ 30$, $T=0.25$ (i.e.
the expiration time  is 3 months=25 years); $r=0.05$ (i.e 5\%
annual), $\s=0.45$ (i.e. 45\% annual). Then $d=d_p=d_c=0.978$,
$\Phi(d)=0.836$, $\mu_c/\mu_p=164/836$. Let  $S_0$ be arbitrary,
$K_p=K_c=S_0$, then $d=d_p=d_c=0.168$, $\Phi(d)=0.567$,
$\mu_c/\mu_p=433/567$.
\section{A consequence for the  seller and a paradox} Consider the
result for a seller (writer) who has sold the options portfolio
described in Theorem \ref{Th1}. The seller has received the
premium $X_0$ and must pay $X_T$ at the time $t=T$.
\par
 Let
 $Y_t$ be the wealth, which was obtained by the seller from   $Y_0=X_0$
by some self-financing strategy. Let $\w Y_T=Y_T-X_T$ be the
terminal wealth after paying obligations to options holder at the
expiration time $T$.
\par
Let us consider  possible actions of the sellers after receiving
the premium. The following strategy  is most commonly presented in
textbooks  devoted to mathematical aspects of option pricing:
\par
{\bf Strategy I}: {\it To replicate the claim $X_T$ using the
replicating strategy.}
\\
As is known, the Black-Scholes price is defined as a minimum
initial wealth such that the option's random claim can be
replicated. For the Strategy I,  the number of shares is
$\Delta(S_t,t)$ at any time $t\in [0,T]$, and $\w Y_T=0$ a.s.,
i.e. {\it there is not neither risk nor any gain}. Thus, it is
doubtful that the seller uses this strategy in practice.
\par
Furthermore,   it was mentioned in Strong (1994)[c. 53] that in
practice the option writer rather just keeps premium as a really
compensation for bearing added risk of for foregoing future price
appreciation or depreciation. Thus, the second strategy is the
following:
\par
{\bf Strategy II}: {\it To invest the premium $X_0$ into bonds,
take no further actions and wait the outcome of the price movement
similarly as the option's holder.}
\par
The seller who sells only put (or only call ) options and uses the
Strategy II puts his/her stake on the random events $K_p\le S_T$
(or $S_T\le K_c$ correspondingly). A gain of the holder implies a
loss for the writer. It looks as a fair game in a case of selling
either put or call separately, because we know that the
Black-Scholes price is fair and chances for gain should be equal
for buyer and seller (otherwise, either  ask or bid will prevail).
But Theorem \ref{Th1} implies that that the seller will receive
non-positive in average gain which is negative in average for any
$a\neq r$ if he/she sells {\it the combination described in
Theorem \ref{Th1}}. In other words, we have a paradox:
\par
\hspace{.5em}
\parbox{.9\textwidth}
{\it A combination of two "fair" deals of selling puts and calls
gives an "unfair" deal.}
\\
\vspace{5mm} 
 The second paradox can be
formulated as following:
\par
  \hspace{.5em}
\parbox{.9\textwidth}
{\it  We know that the Black-Scholes price is
 fair  for buyers
as well as  for sellers (otherwise, either  ask or bid will
prevail). However, we found an options combination such that
buying is preferable, because  the buyer has non-negative and
almost always
 positive average gain,
 but the seller has either  zero  or negative average gain.}
\par
\vspace{5mm}  In practice, unlike as for our generic model,
brokers use sophisticated measures such as insurance, long
positions in stocks or other options, etc. to reduce the risk, but
that does not affect the core of the paradoxes. \par
A possible
explanation is that writer also use the premium $X_0$ to can
receive a gain from $a\neq r$ {\it by using some other strategies}
such as Merton's strategies which are possibly more effective then
options portfolio. In other words, {\it a rational option's seller
does not use neither risk-free replicating of claims nor
"keep-only-bonds" strategy; he/she rather uses strategies which
are able to explore $a\neq r$}.
\section{Proofs} \par
{\it Proof of Theorem  \ref{Th1}.}
  Let $\P_a$ be the
conditional probability measure given $a$. Let $\E_a$ be the
corresponding expectation. We denote $\E_*$ the expectation which
corresponds the risk neutral measure, when $a=r$.
Set
$$
h(a)\defi\E_a X_T.
$$
By the definitions of $X_T$, it follows that
\be
\label{hh}
 h(a)=\mu_p\E_*(K_p-e^{(a-r)T}S_T)^++
\mu_c\E_*(e^{(a-r)T}S_T-K_c)^+. \ee As is known (see e.g. Strong
(1994) and Duffie  (1988)), the Black-Scholes price can be
presented as
\be
\label{13}
\begin{array}{ll}
p_{\BS}(S_0,K,r,T,\s)=e^{-rT}\E_*(K-S_T)^+,
\\
c_{\BS}(S_0,K,r,T,\s)=e^{-rT}\E_*( S_T-K)^+. \end{array}\ee
 Then $$
e^{-rT}h(a)=\mu_pp_\BS(e^{(a-r)T}S_0,K_p,r,T,\s)+ \mu_cc_\BS(
e^{(a-r)T}S_0,K_c,r,T,\s).
$$
By the put and call parity formula, it follows that
$$
e^{-rT}h(a)=\mu_p\biggl[c_\BS( e^{(a-r)T}S_0,K_p,r,T,\s)-
e^{(a-r)T}S_0+K_p\biggr]+ \mu_cc_\BS(e^{(a-r)T}S_0,K_c,r,T,\s). $$
 The following proposition is well known (see e.g. Strong (1994)[p. 100]).
\begin{proposition}
For any $T>0, K>0$, the following holds:
\be
\frac{\p}{\p x}c_{\BS}(x,K,r,T,\s)=\Phi(d),
 \label{(14)} \ee
where $d$ is defined by (\ref{(3)}).
\end{proposition}
Let $d_c(a)$ and $d_p(a)$ be defined as $d_p$ and $d_c$
correspondingly  with substituting $ e^{(a-r)T}S_0$ as $S_0$. Set
$$y=y(a)\defi e^{(a-r)T}, \quad
 R(y)\defi h(T^{-1}\ln y). $$ We have that $$R(y(a))\equiv
h(a),$$
then
$$
 \frac{d R(y)}{d y}=e^{rT}S_0\left[(\mu_p(\Phi(d_p(a))-1)+
\mu_c\Phi(d_c(a))\right]. $$
 By (\ref{(8)}) it follows that
 \be
 \label{(16)} \frac{d R}{dy}(y)|_{y=1}= 0.
\ee
(Note that $y=1$ if and only if $a=1$.) It is know that
$$
\frac{\p ^2}{\p x^2}p_{\BS}(x,K_p,0,T,\s)>0, \quad \frac{\p
^2}{\p x^2}c_{\BS}(x,K_p,0,T,\s)>0
$$
and the derivatives exist. Then $R''(y)>0$ $(\forall y)$, i.e.
$R(\cdot)$ is strongly convex. It follows that $a=r$ is the only
one which minimize $h$. By (\ref{(6)}), (\ref{hh}) it follows
that
$$ R(1)=e^{rT}X_0.
$$
 The uniqueness of the
proportion (\ref{(8)}) follows from the uniqueness of
$\mu_c/\mu_p$ which ensures (\ref{(16)}). This completes the proof
of Theorem \ref{Th1}.
\par
Corrolary  \ref{cor1} follows from (\ref{(9)})--(\ref{(10)}).
Corrolary  \ref{cor2} follows from (\ref{(14)}).

\end{document}